\documentclass[12pt]{article}
\usepackage[cp866]{inputenc}
\usepackage[russian]{babel}
\ifx\pdftexversion\undefined
  \usepackage[dvips]{graphics}
\else
  \usepackage[pdftex]{graphics}
\fi
\usepackage{epsfig}
\usepackage{prep_r1}
\def\Quad{\mathcal{Q}}
\def\cQuad{\overline\Quad}
\def\Cycl{\mathcal{C}}

\def\sign{{\rm sign}\,}
\def\chet{че\-ты\-рех\-уголь\-ник}
\def\mnog{мно\-го\-член}
\def\CAS{{\it Maple}}
\def\fz{{\partial}_z F}
\def\Mr{\mathcal{M}}
\def\Kr{\mathcal{K}}
\def\Dr{\mathcal{D}}
\def\pdq{\partial\Quad}
\def\pd{\partial}
\def\Resl#1#2#3{\mathrm{Resultant}_{#3}(#1,\,#2)}
\newtheorem{Predl}{Предложение}

\begin{document}
\large
\toptitlepage
\thispagestyle{empty}
\vskip 6cm
\large
\centerline{ С.Ю.~Садов}
\vskip 2cm
\centerline{
О НЕОБХОДИМОМ И ДОСТАТОЧНОМ УСЛОВИИ}

\bigskip
\centerline{
ВПИСАННОСТИ ЧЕТЫРЕХУГОЛЬНИКА В ОКРУЖНОСТЬ}
\vskip 10cm
\large
\centerline{\Msk, 2003 г.}
\newpage
\thispagestyle{empty}
\vbox{
\sloppy
\noindent
УДК \quad 514.1+519.11
\newline
\newline
\noindent
 С.Ю.~Садов.
О необходимом и достаточном условии вписанности четы\-рех\-угольника в
окружность.
$\;$
Препринт Института прикладной ма\-те\-ма\-ти\-ки им.~М.В.~Келдыша РАН,
\Msk, 2003 г. \No\ 94.  \newline

Выпуклый четырехугольник со сторонами $a$, $b$, $c$, $d$ и
диагоналями $p$, $q$ является вписанным тогда и только тогда,
когда $\;abp-bcq+cdp-daq=0$. Несмотря на простоту, это условие,
по-видимому, ново
и неожиданно трудно доказуемо. В работе
привлекаются методы компьютерной алгебры и локального нелинейного
анализа.
\vspace{2cm}

\noindent
S.Yu.\ Sadov.
On a necessary and sufficient cyclicity condition for
a quadrilateral.
$\;$
Preprint of the M.V.~Keldysh Institute for Applied Mathematics of
RAS, Moscow, 2003, No\ 94.
\newline

A convex quadrilateral with sides $a$, $b$, $c$, $d$ and diagonals
$p$, $q$ is cyclic iff $\;abp-bcq+cdp-daq=0$.
This condition, in spite of its simplicity, appears to be
unnoted and unexpectedly proof-resilient.
We employ advanced methods of computer algebra and nonlinear analysis.
\newline

\vskip 2cm

Работа выполнена при поддержке
Российского Фонда Фунда\-мен\-тальных Исследований,
гранты 02-01-01067 и 01-01-00517.
\newline

\noindent
\vskip 1cm
E-mail:$\quad$ sadov@keldysh.ru

\vskip 2cm\noindent \IPM{2003}
}

\newpage

\section*{Введение}
В работе рассматривается необходимое и достаточное условие того, что
выпуклый \chet\ $ABCD$ является вписанным (его вершины
лежат на одной окружности). Условие формулируется в терминах длин сторон и
диагоналей \chet{}а.
Обозначим
$$
AB=a,\quad BC=b,\quad CD=c,\quad DA=d,\quad AC=p, \quad BD=q.
$$
Известно классическое условие Птолемея:
\begin{equation}
\label{C2}
 C_2(a,b,c,d,p,q) \eqbydef ac+bd-pq=0.
\end{equation}
(Буква $C$ --- от английского
\ {\it cyclic}\
--- вписанный. Индекс $2$ обозначает степень однородности полинома.)
$\,$
Новое условие задается однородным полиномом третьей степени
\begin{equation}
\label{C3}
C_3(a,b,c,d,p,q) \eqbydef abp-bcq+cdp-daq=0.
\end{equation}
Функция  $C_2$ всегда неотрицательна. В отличие от нее,
знак функции $C_3$ определяет, лежит ли точка $D$ внутри или вне
окружности $ABC$. Такой критерий может быть полезен в
приложениях.

Доказательство достаточности условия (\ref{C3}), найденное
автором, --- не\-эле\-мен\-тар\-но и очень громоздко.  Мы показываем, что
некоторая система поли\-но\-ми\-аль\-ных уравнений не имеет решений в
области допустимых зна\-че\-ний пе\-ре\-мен\-ных.  Существенной частью
доказательства является ло\-каль\-ный анализ системы вблизи границы
области.

Проверка всех выкладок настоящей работы вручную, без ис\-поль\-зо\-ва\-ния
системы компьютерной алгебры
(я ис\-поль\-зо\-вал \CAS), вряд ли возможна.
Более того, вы\-де\-ле\-ние ветвей решения ведет к необходимости
разрешения особенностей и много\-чис\-лен\-ным случаям
и подслучаям. Ана\-лиз некоторых вырожденных слу\-ча\-ев
(см.\ п.~5.1 и 6.1) еще тре\-бу\-ет за\-вер\-ше\-ния и здесь не приводится.

В.П.~Варин 
предложил
значительно более простое доказательство (по-видимому, допускающее
проверку вручную), основанное на открытом им
заме\-ча\-тельном тождестве, связывающем
величины $C_2$ и $C_3$. В оправдание под\-хода, используемого
здесь, укажем на его уни\-вер\-саль\-ность
(сдер\-жи\-ва\-е\-мую недостаточной развитостью программного обеспечения
для ана\-ли\-за особенностей многомерных алгебра\-ических уравнений).

\pagebreak
\vspace*{0.6cm}
\noindent
\kern 20pt
$\lefteqn{\mbox{\includegraphics{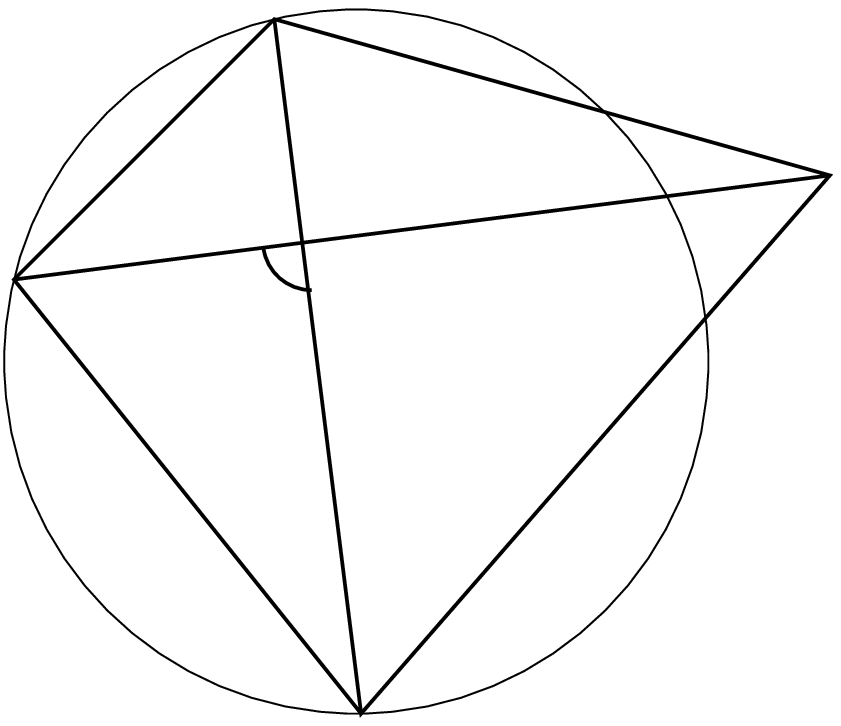}}}{\mbox{\unitlength=1pt
\begin{picture}(437 ,220 )
\put(95 , 0 ){$A$}
\put(-14 ,129 ){$B$}
\put(70 ,216 ){$C$}
\put(242.600 ,159 ){$D$}
\put(45 , 70.100 ){$a$}
\put(32.5 ,175.5 ){$b$}
\put(155.4 ,195.9 ){$c$}
\put(176 ,84 ){$d$}
\put(88 ,178 ){$u$}
\put(100 ,85 ){$x$}
\put(45 ,132 ){$y$}
\put(156 ,146){$z$}
\put(69.5 ,153){$O$}
\put(72 ,121){$\theta$}
\put(115 ,-20){Рис.~1}
\put(300 ,100){$
\begin{array}{c}
p=x+u\\
q=y+z\\
t=\cos\theta
\end{array}
$}
\end{picture}}}
$
\vspace{10pt}

\section{
Критерий Птолемея}

\begin{tm} 
{\rm (А)}$\,$ Для любого \chet{}а $ABCD$ имеет место неравенство
\begin{equation}
\label{Pto}
 ac+bd \geq pq.
\end{equation}
{\rm (Б)}$\,$ Неравенство {\rm(\ref{Pto})} обращается в равенство тогда и
только тогда, когда \chet\ $ABCD$ вписанный.
\end{tm}

Строго говоря, Птолемею принадлежит часть {\it тогда}\
утверждения Б. 

Приведем ссылки на три элементарных доказательства, найденные в легко
доступной литературе:

\smallskip\noindent
1) Доказательство, основанное на неравенстве треугольника для {
пе\-даль\-но\-го} треугольника и условии вырождении этого треугольника в
отрезок {
прямой Симсона} \cite[гл.~2, \S~5-6]{Coxeter}.

\smallskip\noindent
2) Доказательство, основанное на преобразовании инверсии
\cite{Prasol}, задача 28.24.

\smallskip\noindent
3) Доказательство, использующее комплексные числа
\cite{Prasol}, Приложение 1.

Последнее доказательство особенно просто, поэтому позволим себе
вос\-произвести его. Отождествляя векторы на плоскости с комплексными
чис\-ла\-ми и помещая точку $A$ в начало координат, напишем
$$
\vec{AB}=z_1,\quad
\vec{AC}=z_2,\quad
\vec{AD}=z_3,
$$
тогда
$$
\vec{BC}=z_2-z_1,\quad
\vec{CD}=z_3-z_2,\quad
\vec{BD}=z_3-z_1,\quad
$$
и произведения, участвующие в теореме Птолемея, записываются в виде
$$
 ac=
 |z_1 z_3- z_1 z_2|, \qquad
 bd=
 |z_2 z_3- z_1 z_3|, \qquad
 pq=
 |z_2 z_3- z_1 z_2|.
$$
Неравенство (\ref{Pto}) представляет собой неравенство треугольника
для трой\-ки вершин
$\;z_1 z_2$, $\;z_1 z_3$, $\;z_2 z_3$.
Оно обращается в равенство тогда и только тогда, когда точка
$\,z_2 z_3\,$ лежит на отрезке, соединяющем $\,z_1 z_2\,$ и $\,z_1 z_3$.
Это условие выражается формулой
$$
\arg\,(z_1 z_2 - z_2 z_3)\,=\, \arg\,(z_1 z_2 - z_1 z_3).
$$
Преобразуем его к виду
$$
\arg z_2- \arg z_1 \;=\; \arg\,(z_2-z_3)- \arg\,(z_1- z_3)
$$
и перепишем в геометрических обоначениях
$$
 \angle BAC\;=\;\angle BDC.
$$
Получилось условие равенства углов, опирающихся на одну и ту же сто\-ро\-ну
($BC$) нашего \chet{}а, т.е. условие вписанности.
\BlackBox

%
%

\section{
Алгебраическое доказательство теоремы 1 
\\ (Метод грубой силы)}

\subsection{
Пространство выпуклых четырехугольников}

Множество всех шестерок $\,(a, b, c, d, p, q)$, соответствующих выпуклым
не\-вы\-рожденным \chet{}ам, образует пятимерное
под\-мно\-го\-обра\-зие $\Quad$, лежащее в положительном гипероктанте
пространства параметров $\Rr^6$.  Его замыкание $\cQuad$
--- многообразие с кусочно-гладким краем.

Многообразие $\Quad$ допускает удобную глобальную параметризацию.

Обозначим точку пересечения диагоналей $AC$ и $BD$ через
$O$.
Введем пять независимых параметров (см.\ Рис.~1)
\begin{equation}
\label{elem5}
OA=x,\quad\;
OC=u,\quad\;
OB=y,\quad\;
OD=z,\quad\;
\cos \angle AOB=t,
\end{equation}
подчиненных ограничениям
\begin{equation}
\label{rest5}
x >0,\quad\;
u >0,\quad\;
y >0,\quad\;
z >0,\quad\;
-1<t<1.
\end{equation}
Длины сторон и диагоналей даются выражениями
\begin{equation}
\label{elem6}
\begin{array}{l}
a^2=x^2+y^2-2xyt,\\[1ex]
b^2=u^2+y^2+2uyt,\\[1ex]
c^2=u^2+z^2-2uzt,\\[1ex]
d^2=x^2+z^2+2xzt,\\[1ex]
p=x+u\\[1ex]
q=y+z.
\end{array}
\end{equation}
Уравнения (\ref{elem6})
описывают взаимно-однозначное отображение
области (\ref{elem5}) из $\Rr^5$ на
$\,\Quad\subset \Rr_+^6$. Непосредственное описание многообразия $\Quad$
как подмножества в $\Rr_+^6$ в терминах $a,\dots, q$ гораздо более
замысловато. Приведем его для справки. Описание состоит из одного
уравнения и не\-сколь\-ких неравенств
(неравенства треугольников, условия по\-ло\-жи\-тель\-нос\-ти длин
и неравенства, отвечающие за выпуклость).

Четырехугольник есть вырожденный случай тетраэдра, когда все
вер\-ши\-ны лежат
в одной плоскости.
Объем $V$ тетраэдра со сторонами
$\,a$, $b$, $c$, $d$, $p$, $q\,$
дается {\it определителем
Кэли-Менгера} \cite[п.~9.7.3]{Berger}
\begin{equation}
\label{caley}
288\,V=\left|\begin{array}{ccccc}
0 & 1 & 1 & 1 & 1\\
1 & 0 & a^2 & p^2 & d^2 \\
1 & a^2 & 0 & b^2 & q^2 \\
1 & p^2 & b^2 & 0 & c^2 \\
1 & d^2 & q^2 & c^2 & 0
\end{array}
\right|
\end{equation}
(Об истории формулы объема тетраэдра см.~\cite{Sabitov}.)
Раскрывая определитель, получим уравнение \chet{}а в явном виде
\begin{equation}
\label{quadeq}
\begin{array}{l}
a^4 c^2 + a^2 c^4 + b^4 d^2  + b^2 d^4 + p^4 q^2 + p^2 q^4 +
\\[2ex]
+ (abp)^2 + (bcq)^2 +(cdp)^2 +(daq)^2 -
\\[2ex]
-(abc)^2 - (abd)^2 - (acd)^2 - (acp)^2 - (acq)^2 - (apq)^2 -
\\[2ex]
-(bcd)^2  - (bdp)^2- (bdq)^2 - (bpq)^2 - (cpq)^2 - (dpq)^2
\\[2ex]=0.
\end{array}
\end{equation}
К этому уравнению добавляются упомянутые неравенства,
часть из ко\-то\-рых очевидна
\begin{equation}
\label{quadineq}
\begin{array}{l}
a> 0,\;\quad
b> 0,\;\quad
c> 0,\;\quad
d> 0,\;\quad
p> 0,\;\quad
q> 0,
\\[0.5ex]
|a-b|<p<a+b,\qquad
|c-d|<p<c+d,
\\[0.5ex]
|b-c|<q<b+c,\qquad
|a-d|<q<a+d.\qquad
\end{array}
\end{equation}
Неравенства, отвечающие за выпуклость, не столь тривиальны.
При за\-дан\-ных $a$, $b$, $c$, $d$, $p$ построим треугольники
$ABC$ и $ADC$ на общем ос\-но\-ва\-нии $AC$. Вершины
$B$ и $D$ могут лежать либо по
разные стороны от прямой $AC$,
либо по одну сторону (Рис.~2). Найдем $q=BD$ из построенной
конфигурации.

\bigskip\bigskip
$\lefteqn{\mbox{\includegraphics{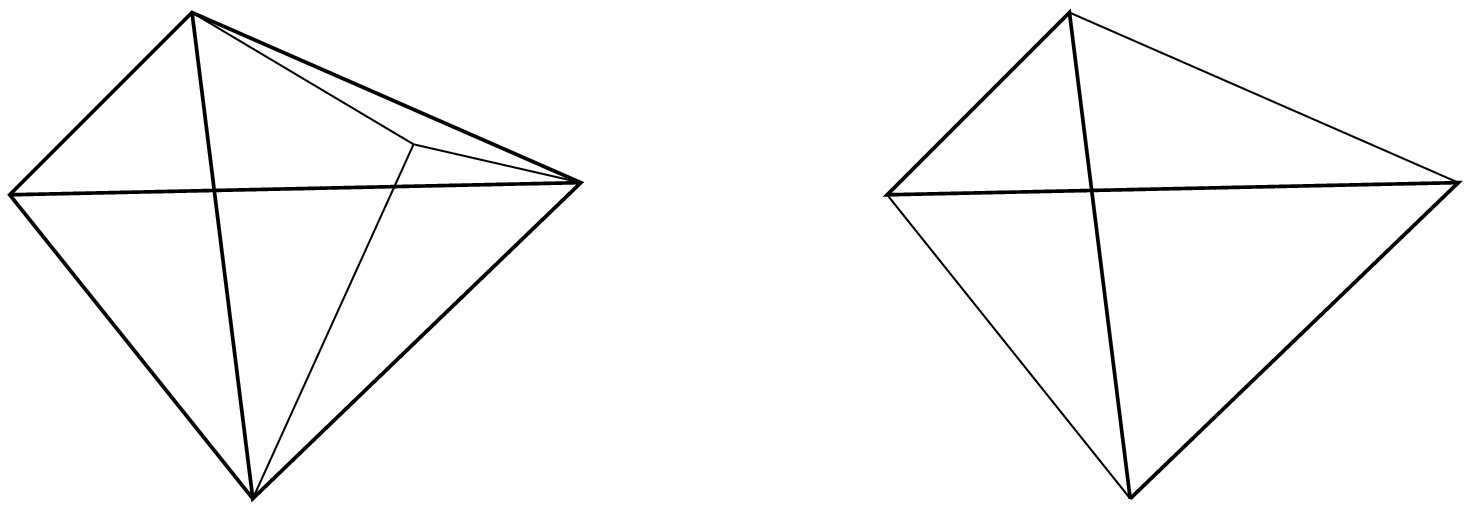}}}{\mbox{\unitlength=0.7pt
\begin{picture}(630.000 ,220.000 )
\put(88, -2 ){$A$}
\put(-4, 148 ){$B$}
\put(143 ,145 ){$B'$}
\put(80 ,218 ){$C$}
\put(239 ,120 ){$D$}
\put(115 ,-28){Рис.~2}
\put(450 , -2 ){$A$}
\put(358, 148 ){$C$}
\put(441 ,218 ){$B$}
\put(600 ,120 ){$D$}
\put(476 ,-28){Рис.~3}
\end{picture}}}
$

\vspace{26pt}
\noindent
В  обоих случаях условия (\ref{quadeq}) и (\ref{quadineq})
выполнены.
Неравенство, отли\-ча\-ющее  выпуклый \chet\ $ABCD$ от невыпуклого
$AB'CD\,$:
$$
\angle BAD\; > \;\angle CAD.
$$
Используя теорему косинусов, перепишем это условие в виде
$$
\frac{a^2+d^2-q^2}{ad}\;=\;2\cos\angle BAD \;<\; 2\cos\angle CBD
\;=\;\frac{a^2+p^2-b^2}{ap}.
$$
Приводя дроби к общему знаменателю, получим полиномиальное
не\-ра\-вен\-ство. Оно еще не гарантирует выпуклости \chet{}а
$\,ABCD$. Необходимо добавить неравенство, исключающее случай,
когда точки $A$ и $C$ лежат по одну сторону от прямой $BD$.
Двух неравенств все еще недостаточно. Например,
система неравенств
$$
 \angle BAC < \angle BAD, \qquad
 \angle ABD < \angle ABC
$$
допускает конфигурацию (четырехвершинник) с неправильным порядком вершин
(Рис.~3). Трех неравенств
$$
 \angle BAC < \angle BAD, \qquad
 \angle ABD < \angle ABC, \qquad
 \angle CBD < \angle ABC
$$
уже достаточно для характеризации выпуклых \chet{}ов
т.е.\ к (\ref{quadineq}) добавляются неравенства
\begin{equation}
\label{conv}
\begin{array}{l}
p\,(a^2+d^2-q^2)\,<\,  d\,(a^2+p^2-b^2), \\[1.5ex]
q\,(a^2+b^2-p^2)\,<\,  a\,(b^2+q^2-c^2), \\[1.5ex]
q\,(a^2+b^2-p^2)\,<\,  b\,(a^2+q^2-d^2).
\end{array}
\end{equation}
Можно ли упростить условия выпуклости --- уменьшить число неравенств
или понизить их степень (даже увеличив их количество), принимая во внимание
(\ref{quadineq}), --- я не знаю.

\subsection{
Пространство вписанных четырехугольников}

Шестерки длин сторон и диагоналей {\it вписанных}\ \chet{}ов
обра\-зуют подмногообразие $\Cycl \subset \Quad$ коразмерности 1.
В координатах $x,\dots,t$ оно описывается простым уравнением
второй степени
\begin{equation}
\label{Cycl0}
Cycl(x,y,z,u,t)\, \eqbydef \,xu-yz=0.
\end{equation}
(Элементарная теорема: треугольники $AOB$ и $COD$ подобны по двум
углам $\,\angle AOB=\angle COD\,$ и $\,\angle OAB=\angle ODC\,$ тогда
и только
тогда, когда точки $A$, $B$, $C$, $D$ лежат на одной окружности.)

С другой стороны, часть Б теоремы Птолемея утверждает, что
в коор\-ди\-на\-тах $\,a,\dots,q$ подмногообразие $\Cycl$
выделяеется в $\Quad$ условием (\ref{C2}). Цель этого пункта ---
установить эквивалентность (\ref{C2}) и (\ref{Cycl0}) алгебраически.
За\-одно мы докажем и неравенство (\ref{Pto}).

Выразив $a,\dots,q$ из (\ref{elem6}) и подставив в (\ref{C2}), получим
выражение,
содер\-жа\-щее квадратные корни. Домножая на сопряженные выражения,
изба\-вим\-ся от иррациональностей и получим
$$
\begin{array}{l}
 C_2\cdot (ac-bd-pq)(ac+bd+pq)(ac-bd+pq)=
 \\[2ex]
 - 4u^4 x^2 z^2  + 4x^4 u^2 z^2 t^2  + 32x^2 y^2 u^2 z^2  + 4y^4 u^2 z^2 t^2
 + 4x^2 y^2 t^2 u^4  + 4 x^2 y^2 t^2 z^4  \\[1ex]
 - 4x^2 y^2 z^4  - 4y^4 u^2 z^2  - 4x^2 y^2 u^4
 - 4x^4 u^2 z^2  - 4y^4 x^2 z^2
- 4u^2 y^2 z^4  - 4u^2 y^2 x^4  \\[1ex]
 - 32x^2 y^2 u^2 z^2 t^2  + 16x^3 yt^2 u^3 z
 - 8x^3 yt^2 z^3 u - 8y^3 xt^2
u^3 z + 16y^3 xt^2 z^3 u \\[1ex]
 + 4u^4 x^2 z^2 t^2  + 4y^4 x^2 z^2
t^2  + 4u^2 y^2 t^2 x^4  + 4u^2 y^2 t^2
z^4  + 8x^3 uy^3 z + 16x^3 uy^2 z^2  \\[1ex]
 + 8x^3 uyz^3  + 16x^2 u^2 y^3 z + 16
x^2 u^2 yz^3  + 8xu^3 y^3 z + 16xu^3
y^2 z^2  + 8xu^3 yz^3  \\[1ex]
 + 8x^3 u^3 z^2 t^2  + 8y^3 u^2 z^3 t^2  + 8x^3 y^2 t^2 u^3
 + 8x^2 y^3 t^2 z^3  - 8x^2 u^4 yz - 8x^4 u^2 yz \\[1ex]
 - 16x^3 u^3 yz - 8y^4 z^2 xu - 16 y^3 z^3 xu - 8y^2 z^4 xu -
 8u^3 yt^2 z^3 x - 8y^3 ut^2 x^3 z \\[1ex]
 - 8x^3 u^3 y^2  - 8x^3 u^3 z^2  - 8x^2 z^3 y^3  - 8y^3 z^3 u^2
 - 16y^2 u^3 z^2 t^2 x + 8y^4 uz^2 t^2 x \\[1ex]
 - 16y^3 u^2 zt^2 x^2  + 8x^2 y t^2 u^4 z
- 16x^3 uz^2 t^2 y^2  + 8x^4 u^2 z t^2 y - 16x^2 u^2 z^3 t^2 y \\[1ex]
 + 8xy^2 t^2 z^4 u .
\end{array}
$$
В дальнейшем громоздкие полиномы,
возникающие в вычислениях, не выписываются явно.
Те из них, которые используются
лишь локально, в
пределах конкретного этапа рас\-суж\-дения,
обо\-зна\-ча\-ют\-ся $\,P_k$, где $\,k\,$ --- общая степень.
Так, назовем приведенный выше мно\-го\-член $P_{10}$.
(Его
\linebreak
содержательная характеризация
как главного минора опре\-де\-ли\-теля Кэли-Менгера
(\ref{caley})$\,$
\cite[п.~9.7.3.8]{Berger} здесь не используется.)

Компьютерная факторизация приводит к простому
результату, ко\-то\-рый запишем в смешанных переменных,
связанных соотношениями (\ref{elem6}),
\begin{equation}
\label{facptol}
 P_{10}(x,y,x,u,t)\,=-\,4(1-t^2) p^2 q^2 (xu-yz)^2.
\end{equation}
Поскольку
\begin{equation}
\label{facptin}
(1-t^2) p^2 q^2 \,>0,
\end{equation}
то (\ref{C2}) влечет (\ref{Cycl0}).
Доказано утверждение {\it только тогда}\ 
те\-оре\-мы~1 (Б).

Предположим, что неравенство (\ref{Pto}) неверно для некоторого \chet{}а.
Тогда существует кривая $\gamma:\, [0,1]\to \Quad$ в пространстве \chet{}ов,
такая, что $C_2(\gamma(0))>0$ и $C_2(\gamma(1))<0$,
следовательно, существует точка $s\in (0,1)$, в которой
$C_2(\gamma(s))=0$.
В силу (\ref{facptol}) и (\ref{facptin}), для любого невырожденного
\chet{}а
\begin{equation}
\label{factor2}
(ac+bd-pq)(ac+pq-bd)(bd+pq-ac)\geq 0.
\end{equation}
Если 2-й и 3-й сомножители не обращаются в 0
в $s_0$, то произведение (\ref{factor2}) меняет знак, что
невозможно. Следовательно, хотя бы два из со\-мно\-жи\-те\-лей в
(\ref{factor2}) одновременно равны 0 в $s_0$. Однако тогда
\chet\  $\gamma(s_0)$ оказывается вырожденным.
Действительно, если, например
$$
ac+bd-pq=ac+pq-bd=0,
$$
то $ac=0$.  Полученное противоречие доказывает справедливость
не\-ра\-вен\-ства (\ref{Pto}).

Остается доказать утверждение {\it тогда}\ части Б. Для этого
следует уста\-но\-вить, что $(ac+pq-bd)$ и $(bd+pq-ac)$ не могут
(по отдельности) обра\-щать\-ся в $0$.
Многообразие $\Cycl$ линейно связно:
любой вписанный \chet\ можно деформировать
в любой другой, вписанный в ту же окруж\-ность,
сдвигом вершин вдоль окружности, при котором вершины
никогда не сли\-ва\-ются.
В каждой точке многообразия $\Cycl$ обращается в 0 ровно один
из множителей в (\ref{factor2}) --- иначе \chet\
был бы вырожденным. Множество нулей каждого
из них замкнуто в $\Cycl$.
Из связ\-но\-сти следует, что два из этих множеств
пусты, а оставшееся со\-в\-па\-да\-ет с $\,\Cycl$.
Ясно (из любого примера вписанного \chet{}а), что
не пусто множество нулей функции $C_2$, значит,
$\,(Cycl=0)\,\Leftrightarrow\,(C_2=0)$.
\BlackBox

\section{
Кубический критерий вписанности}

\begin{tm}
Выпуклый \chet\ $\,ABCD\,$ является вписанным
тогда и только тогда, когда выполнено условие\ {\rm(\ref{C3})}.
Если $\,C_3\neq 0$, то
\begin{equation}
\label{C3sign}
\sign\, C_3=\sign\, Cycl.
\end{equation}
\end{tm}

\noindent
{\bf Замечания}

\smallskip
\noindent
1. {\it Существенность условия выпуклости.}\
В отличие от условия Птолемея (\ref{C2}), условие (\ref{C3})
выполняется для некоторых невыпуклых (и, значит, невписанных)
\chet{}ов. В Приложении вычислено частное \linebreak
одно\-параметрическое семейство таких \chet{}ов.

\smallskip\noindent
2. {\it Знак\ $C_3\,$}\ определяет положение вершины \chet{}а
от\-но\-си\-тель\-но
окружности, проведенной через три другие вершины.
Вероятно, $C_3$ --- простейшая
функция параметров $a,\dots,q$, обладающая этим свой\-ством.
Cледующие условия равносильны:
$$
\label{C3ineq}
\begin{array}{rl}
 {\rm(i)}\quad\; &
\mbox{\rm точка $A$ лежит внe окружности $BCD$}
\\[0.5ex]
 {\rm(ii)}\quad\; &
\mbox{\rm точка $C$ лежит внe окружности $BAD$}
\\[0.5ex]
 {\rm(iii)}\quad\; &
\mbox{\rm точка $B$ лежит внутри окружности $ADC$}
\\[0.5ex]
 {\rm(iv)}\quad\; &
\mbox{\rm точка $D$ лежит внутри окружности $ABC$}
\\[0.5ex]
 {\rm(v)}\quad\; &
 \sgn\,(xu-yz)\,=\sgn\,C_3\,>0;
\end{array}
$$

\medskip\noindent
{\bf Доказательство необходимости}$\;$ условия (\ref{C3}) элементарно.
Пусть \chet\ $ABCD$ вписан в окружность радиуса $R$.
Тогда
$$
 4R=\frac{abp}{S_{ABC}}=\frac{bcq}{S_{BCD}}=\frac{cdp}{S_{CDA}}=
 \frac{daq}{S_{DAB}},
$$
где $S$ --- площадь треугольника. Равенство (\ref{C3}) следует
выражает равно\-со\-став\-лен\-ность нашего \chet{}а:
$$
 S_{ABCD}=S_{ABC}+S_{BCD}
 =S_{CDA}+S_{DAB}.
\eqno\mbox{\BlackBox}
$$

\medskip\noindent
{\bf Доказательство
достаточности}$\;$ гораздо сложнее. Метод грубой силы в простой
редакции недостаточен: при рационализации условия $C_3=0$ появляются
паразитные решения --- см.\ п.~5.2.
Изложим план до\-ка\-за\-тель\-ства.
Рас\-смот\-рим величину $C_3$ как функцию переменной $z$ при
фикси\-ро\-ван\-ных
$x,u,y,t$. Пишем
$$
 C_3(a,b,c,d,p,q)=F(x,y,z,u,t)
 \;\quad\mbox{с подстановкой (\ref{elem6})}.
$$
Тогда

1. $\quad F\to-\infty\;\;$ при $z\to \infty$. (Простой факт, Лемма 4.1
ниже.)

2. $\quad F>0\;\;$ при $z=0$. (Предложение 4.2, трудное.)

3. $\quad F=0$, когда точка $D$, двигаясь по лучу OD (меняем
$z$), по\-па\-да\-ет на окружность ABC. При этом
$\,z=z_{*}=xu/y$.

Остается показать, что функция $F(\dots,z)$ не имеет других нулей
на $(0,\infty)$.
Это легко установить в частном случае, осесимметричных \chet{}ов с
перпендикулярными диагоналями (Лемма 3.1). Таким об\-ра\-зом,
существуют значения $x$, $y$, $u$, $t$, при которых функция
$F(\dots, z)$ имеет единственный положительный корень.
Смена числа корней может про\-изойти только с появлением кратного корня,
т.е.\ решения системы
\begin{equation}
\label{fsys}
F(x,y,z,u,t)=0, \qquad
\fz(x,y,z,u,t)=0.
\end{equation}

Предположим, что множество решений системы
(\ref{fsys}) на $\Quad$ непусто. Обозначим его
$\Mr$. Изучая асимптотики решений
системы (\ref{fsys}) вблизи границы
$\pdq$ множества
$\Quad$, мы найдем, что не существует семейства
ре\-ше\-ний в $\Quad$, имеющих предельную точку
на $\pdq$. (Хотя на самой границе есть много решений.)
$\;$ Этот анализ (с пробелом, указанном во Введении)
изложен в заключительном \S\,6.

Сделаем редукцию по однородности и рассмотрим
сечение $\Quad_1$, вы\-де\-лен\-ное
в $\Quad$ уравнением
\begin{equation}
\label{Quad1}
 x+y+z+u=1.
\end{equation}
Множество $\Mr_1=\Mr\cap \Quad_1$
компактно, поскольку его замыкание не пере\-се\-ка\-ется
с $\pdq_1$.
Максимум функции $f(x,y,z,u,y)=t$ на $\Mr_1$ достигается.
Точка максимума находится методом множителей Лагранжа,
который
\linebreak
при\-во\-дит (Лемма 5.1) к усиленной системе по сравнению
с (\ref{fsys}): либо
\begin{equation}
\label{gradsys}
F\,=\,F_x\,=\,F_y\,=\,F_z\,=\,F_u\,=\,0,
\end{equation}
либо $\,F=F_z=F_{zz}=0$.
В \S~5 мы докажем, что система (\ref{gradsys}) не имеет
решений в области (\ref{rest5}). Отсутствие решения во втором случае пока
под\-твер\-жде\-но лишь численным сканированием.

Таким образом, множество решений системы
(\ref{fsys}) в области (\ref{rest5}) пусто,
следовательно, при любых значениях
$x$, $y$, $u$, $t$ корень $z=z_*$ функции
$F(\dots,z)$
--- единственный.
\BlackBox

\begin{Lemma}
В случае $\,x=u\,$ и $\,t=0\,$ имеем $\,\fz<0\,$ при $\,z>0$.
\end{Lemma}

\Proof $\;$
В данном случае $\,a=b\,$, $\;c=d=\sqrt{x^2+z^2}$,
$\;\pd_z c\,=z/c$, 
$$
\begin{array}{l}
  F\,=\,2x\,(a^2+c^2)\,-\,2ac\,(y+z),\\[1.5ex]
  \fz\,=\,4xz\,-\,2ac\,-\,2a\,({z}/{c})\,(y+z)
  \,<\,2a\,(2z- (c^2+z^2)/c)\,\leq 0.
\qquad\mbox{\msam \symbol{4}}
\end{array}
$$

\rem{ Denote $AD=1$, $AB=z_1$, $AC=z_2$, $v_1=z_1 z_2(z_1-z_2)$,
$v_2=z_2(z_2-1)$, $v_3=z_1(z_1-1)$, $v_4=(z_1-z_2)(z_2-1)(z_1-1)$.
Then $v_1+v_2=v_3+v_4=z_2(z_1-1)(z_1-z_2+1)$,
$$
 C_3=|v_1|+|v_2|-|v_3|-|v_4|.
$$
}

\section{
Поведение функции $F$ при $z\to\infty$ и при $z\to0^+$}

\subsection{Результаты}
В настоящей работе относительно простые утверждения называются
Лем\-ма\-ми, более сложные --- Предложениями.

\begin{Lemma}
Пусть $\,x$, $y$, $u$, $t\,$ фиксированы. Тогда при $z\to\infty$
$$
 F(x,y,u,z,t) = (p-a-b)z^2 + O(z).
$$
В частности, $F<0$ при больших $z$.
\end{Lemma}

\Proof $\;$
Поскольку $\,c=z+O(1)$, $\,d=z+O(1)$, $\,q=z+O(1)$, имеем
$$
F(x,u,y,z,t)\,=\, 
(ab+cd)p-(ad+bc)q\,=\,z^2 p-(az+bz)z+O(z).
$$
\BlackBox

\begin{Predl}
Пусть точка $D$ лежит на стороне $AC$ тре\-уголь\-ни\-ка $ABC$
и отлична от точек $A$ и $C$.
Тогда для вырожденного \chet{}а $ABCD$ имеем
$C_3(a,b,c,d,p,q)>0$ при отсутствии прочих вырождений.
В более явном виде,
$$
 (AB\cdot BC+AD\cdot CD)\cdot AC >
 (AB\cdot AD+CB\cdot CD)\cdot BD.
$$
\end{Predl}

\subsection{
План доказательства Предложения 4.2}

Не ограничивая общности, положим $p=1$ и примем $\angle ADB\leq \pi/2$.
На множестве $\;\Dr=\{0<x<1,\;y>0,\; 0\leq t<1\}\;$
рассмотрим функцию
$$
 F_0(x,y,t)\,\eqbydef\, F(x,y,\,z=0,\,u=1-x,\,t).
$$
Множество $\Dr$ связно. Достаточно доказать, что на нем  $\,F_0\neq 0$:
тогда знак $F_0$ постоянен, а из Леммы 3.1 следует, что $F_0>0$
при $t=0$, $x=1/2$. В п.~4.3 рассматривая рационализацию условия
(\ref{C3}), покажем, что $F_0\neq 0$ при $y\geq 1$.
Затем в п.~4.4 исследуем $F_0$
на границе компактного множества
$\Kr$, заданного неравенствами
$$
0\leq x\leq 1, \qquad 0\leq t\leq 1, \qquad 0\leq y\leq 1,
$$
и покажем, что там $F_0\geq 0$. В п.~4.5 установим, что
уравнение $\nabla F_0$ не имеет решений во внутренности
множества $\Kr$, следовательно, ${\rm min\,}F_0$ достигается на его границе
и $F_0>0$ на внутренности $\Kr$.
\BlackBox

\subsection{Рационализация и
неравенство $\,F_0[y\geq 1]\neq 0$}

Выражение $F_0$ через $a,\dots,q$
содержит квадратичные иррациональности.
Функция
$$
\begin{array}{l}
 R(a,b,c,d,p,q)\,=\,C_3(a,b,c,d,p,q)\, C_3(a,b,-c,-d,p,q)\,
 \\[0.5ex]
 \times\,C_3(-a,b,c,d,p,q)\, C_3(-a,b,-c,-d,p,q)
\end{array}
$$
есть рациональная функция от $a^2$, $b^2$, $c^2$, $d^2$, $p$, $q$,
и, следовательно, ра\-ци\-о\-наль\-ная функция от $\,x,u,y,z,t$.
Ее ограничение на рассматриваемое множество
$z=0$, $p=1$, $a=x$, $b=u=1-x$
обозначим через $R_0$.

С помощью \CAS\ находим факторизацию
$$
 R_0=4\,(1-t^2)\,y^2\, x \, T_0(x,y,t)
$$
где $T_0$ --- многочлен.
При $|t|\neq 1$, $\;x,y\neq 0\;$ неравенство $\, T\neq 0\,$ влечет
$\,R_0\neq 0\,$ и, следовательно, $\,F_0\neq 0$.
Многочлен $T_0$ представим в виде
$$
\begin{array}{l}
 T_0(x,y,t)\,=\, A(y^2-1)+4y^2(1-t^2),
 \\[1.5ex]
 A=3y^2+ 4(1-2x)ty + (1-2x)^2.
 \end{array}
$$
Поскольку при $|t|< 1$,
$$
A\,>\, 3y^2- 4y |1-2x| + (1-2x)^2\;=(3y-|1-2x|)(y-|1-2x|)
$$
то при $\,y\geq 1$, $\;0<x<1$, имеем $\;A> 0$, откуда $\,T_0>0$.

\subsection{
Предельные случаи $\;x=0,\,1$, $\;z=0\;$ и $\;t=1$}

\noindent
а)  При $x=0$ (случай $x=1$ аналогичен) имеем
$\,d=0$, $\;c=1$, $\;a=q$, $\;\,$ $F_0=abp-bcq=ab-ba=0$.

\smallskip\noindent
б) При $z=0$
имеем $\,q=0$, следовательно $\;F_0=ab+cd>0$.

\smallskip\noindent
в) При $t=1$ имеет место результат более общего характера
(без пред\-по\-ло\-же\-ния $z=0$), который неоднократно понадобится в
дальнейшем.  Из него следует, что $F_0|_{t=1}\geq 0$.

\begin{Lemma} Имеет место формула
\begin{equation}
\label{C3t1}
 F|_{t=1}\;=(x+y+z+u)(x-y)(u-z)\,(\sgn(x-y)\,+\,\sgn(u-z)).
\end{equation}
\end{Lemma}

\Proof$\;$ Формула (\ref{C3t1}) следует из равенств при $t=1$:
$$
 a=|y-x|,\quad\;
 b=y+u,\quad\;
 c= |u-z|,\quad\;
 d= z+x.
\eqno\mbox{\BlackBox}
$$

\subsection{
Отсутствие внутренних точек экстремума функции $F_0$}

В этом пункте мы используем декартовы координаты.
Поместим начало ко\-ор\-ди\-нат в точку $A$ и направим первую ось вдоль
$AC$.  Ко\-ор\-ди\-на\-ты остальных точек суть $\,D(d,0)$, $C(p,0)$,
$B(\xi,\eta)$, $\;\eta>0$.
Параметры $\;d$, $\,c$, $\,p=c+d\;$ не зависят от $\xi, \eta$,
а параметры $a$, $b$, $q$ и их производные даются формулами
$$
\begin{array}{lclcl}
\dst a^2=\xi^2 +\eta^2 &\qquad&
\dst a_\xi=\frac{\xi}{a} &\qquad&
\dst a_\eta=\frac{\eta}{a}
\\[2ex]
\dst b^2=(\xi-p)^2 +\eta^2 &\qquad&
\dst b_\xi=\frac{\xi-p}b &\qquad&
\dst b_\eta=\frac{\eta}b
\\[2ex]
\dst q^2=(\xi-d)^2 +\eta^2 &\qquad&
\dst q_\xi=\frac{\xi-d}{q} &\qquad&
\dst q_\eta=\frac{\eta}q
\end{array}
$$
Отсюда
\begin{equation}
\label{gradsysz0}
\begin{array}{l}
\dst
\pd_{\xi} F_0\,=\,
\frac{(bp-dq)\xi}{a}+
\frac{(ap-cq)(\xi-p)}{b}-\frac{(ad+bc)(\xi-d)}{q},
\\[2ex]                          `
\dst
\pd_{\eta}F_0\,=\,
\eta\left(\frac{(bp-dq)}{a}+
\frac{(ap-cq)}{b}-\frac{(ad+bc)}{q}\right).
\end{array}
\end{equation}
Назовем набор величин $(a,\dots,q,\xi,\eta)$ допустимым, если
он соответствует конфигурации, описанной в Предложении 4.2
и треугольник $ABC$ не\-вы\-рож\-ден.

\begin{Lemma} 
Система $\;\pd_{\xi} F_0=\pd_{\eta} F_0=0\;$ не имеет
допустимых решений.
\end{Lemma}

\Proof$\;$
Вычитая из первого уравнения в (\ref{gradsysz0}) второе, умно\-жен\-ное
соответственно на $\,\xi/\eta\,$ и $\,(\xi-p)/\eta$,
$\,$ учитывая, что $\,c+d=p$,
$\,$ и из\-бав\-ля\-ясь от знаменателей,
упростим заданные уравнения:
\begin{equation}
\label{simgrad0}
\begin{array}{l} \dst
pq\,(cq-ap)\,+\,bd\,(ad+bc)\,=\,0,
\\[1ex]
\dst
pq\,(dq-bp)\,+\,ac\,(ad+bc)\,=\,0.
\end{array}
\end{equation}
Покажем, что допустимые решения удовлетворяют равенству $a=b=q$.

Беря сумму уравнений и подставляя $\,p=c+d$,
получим факторизуемое уравнение:
$$
[(b-q)d\,+\,(a-q)c]\,[(b-q)c\,+\,(a-q)d]\,=0.
$$
Возможны два случая.

Случай I: $\;(b-q)d+(a-q)c=0$.
В качестве второго уравнения возьмем
определитель системы (\ref{simgrad0}) относительно неизвестных
$pq$ и $(ad+bc)$,
$$
0=\Delta\,=\,ac(cq-ap)+bd(bp-dq)\,=
\,ac\,[(q-a)c-ad]\,+\,bd\, [(b-q)d+bc].
$$
В данном случае
$$
 \Delta\,=\,
\,ac\,[(b-q)d-ad]\,+\,bd\, [(q-a)c+bc]=cd(b-a)(a+b+q),
$$
откуда $a=b$ на допустимом решении. Далее,
$\;0=(b-q)d+(a-q)c=(b-q)(c+d)$, следовательно, $q=b=a$.

Случай II: $\;(b-q)c+(a-q)d=0$, т.е.\ $\,ad+bc=pq$.
Подставляя в первое уравнение
системы (\ref{simgrad0}), получим
$
\;pq\,(cq-ap+bd)=0.
$
С учетом $p=c+d$, получилась однородная линейная система уравнений
$$
 (b-q)c+(a-q)d=0,\qquad
 (q-a)c+(b-a)d=0
$$
относительно $c$ и $d$ с определителем
$$
 (b-q)(b-a)+(a-q)^2= (a-q)^2-(a-q)(b-q)+(b-q)^2.
$$
Равенство нулю возможно лишь при $a=b=q$.

Итак, в обоих случаях (\ref{simgrad0}) влечет $\,a=b=q$.
Однако равенство $a=b=q$ геометрически невозможно,
т.к. всегда $\;q<\,{\rm max\,}(a,b)$.
Этим завершается доказательство Леммы 4.4 и Предложения 4.2.
\BlackBox

\section{Несуществование компактной компоненты \\ подмножества
$\,F=F_z=0\,$ в $\,\Quad_1$}

\subsection{Усиление системы (\ref{fsys})}

\begin{Lemma}
Предположим, что множество решений системы {\rm(\ref{fsys})}
на множестве $\Quad_1\,$ {\rm (см.\ (\ref{rest5}), (\ref{Quad1}) )}
компактно (т.е.\ не имеет предельных точек на $\pdq_1$) и непусто.  Тогда
хотя бы одна из систем {\rm(\ref{gradsys})} или
\begin{equation}
\label{Fzz}
 F\,=\,F_z\,=\,F_{zz}=0
\end{equation}
 имеет решение в $\Quad$.
\end{Lemma}

\Proof $\;$
Из предположений Леммы следует, что множество $\tilde\Mr_1$ решений
системы (\ref{fsys}) на множестве $\,\tilde\Quad_1=\{u=1\}\cap\Quad\,$
также компактно и непусто. Можно считать, что $F$ теперь
является функцией лишь че\-ты\-рех переменных $x$, $y$, $z$, $t$.
Рассмотрим оптимизационную задачу
$$
 t|_{\tilde\Mr_1}\,\to\,\max.
$$
Составим выражение Лагранжа
\begin{equation}
\label{LaGr}
 t-\lambda F - \mu  F_z.
 \end{equation}
Дифференцируя по $z$ и учитывая (\ref{fsys}), получим
$\mu F_{zz}=0$.

В случае $F_{zz}=0$ имеем систему (\ref{Fzz}). Если же
$F_{zz}\neq 0$, то $\,\mu=0$. $\,$ Дифференцируя
(\ref{LaGr}) по $t$, находим $1-\lambda F_t=0$, откуда
$\lambda\neq0$.  Дифференцируя (\ref{LaGr}) по $x$ и $y$,
находим теперь $F_x=F_y=0$. Мы получили 4 из 5 уравнений
(\ref{gradsys}). Оставшееся уравнение $F_u=0$ следует из
первых четырех в силу тождества Эйлера для однородных функций.
\BlackBox

В оставшейся части параграфа мы докажем несуществование
решений системы (\ref{gradsys}), оставляя за рамками данной
работы анализ системы (\ref{Fzz}). Этот пробел упомянут во
Введении.

\subsection{Рационализация}

Функция $\,F(x,\dots,t)$, являющаяся
композицией $\,C_3(a,\dots,q)\,$ с под\-ста\-нов\-ка\-ми (\ref{elem6}),
содержит квадратные корни.

\medskip\noindent\underline{Частичная рационализация}.$\;$
Положим $\,F^{*}=(ab+cd)p+(ad+bc)q\,$ и введем выражение
\begin{equation}
\label{ratR}
\begin{array}{c}
 R(x,y,z,u,t)= F\cdot F^{*}= R_0(x,\dots,t)
  +\lambda R_1(x,\dots,t),
\\[1ex]
 \lambda=abcd.
\end{array}
\end{equation}
Здесь $\;R_1=2(p^2-q^2)$,
$\;R_0=(a^2 b^2+c^2 d^2)p^2-(a^2 d^2+b^2 c^2)q^2\;$ ---
многочлен степени 2 по $\,t\,$ и однородный степени 6 по переменным
$\,x,\dots,u$.
\linebreak
Его разложение содержит 64 члена.
Очевидно, что на множестве $\Quad$ урав\-не\-ния $F=0$ и $R=0$
экви\-валентны, поскольку $\,F^{*}>0$.

\medskip\noindent\underline{Полная рационализация}.$\;$
Функция
$$
 RR(x,y,z,u,t)=(R_0+\lambda R_1)(R_0 -\lambda R_1)
$$
является многочленом по переменным $x,\dots,t$. Имеем факторизацию
\begin{equation}
\label{facR}
RR\,=\, 4\,(t^2-1)\,p^2\, q^2\;Cycl(x,y,z,u)\;T(x,y,z,u,t),
\end{equation}
где $T$ --- многочлен с 62 членами, квадратичный по $t$ и
однородный степени 6 по остальным переменным.
Приводим этот многочлен в явном виде,
следуя принципу "Лучше один раз увидеть\ldots" :
$$
\begin{array}{ll}
T\,=& 4 T_2\, t^2\,+\,4 T_1\,t\,+ T_0,\\[2ex]
T_0\,=&
2y^3z^3+3yzx^4+2y^3zx^2+2yz^3x^2-2x^3uy^2-2xu^3y^2-3xuy^4-
\\[1ex] &
2xu^3z^2+2x^2yzu^2-zy^5-yz^5+ux^5-2xuy^2z^2-3xuz^4-2x^3uz^2
\\[1ex] &
+xu^5+3yzu^4+2y^3zu^2+2yz^3u^2-2x^3u^3,
\end{array}
$$
$$
\begin{array}{ll}
T_1=&
z^2y^3x-y^2x^3z+y^4xz-y^2z^3x+yz^2x^3-z^4xy-x^3yu^2+x^2y^3u
\\[1ex] &
-x^4uy+x^2u^3y-xu^2y^3+z^4uy+y^2z^3u-yz^2u^3-y^3uz^2-y^4uz
\\[1ex] &
+u^3y^2z-u^4xz+xu^2z^3-x^2u^3z-x^2uz^3+x^4uz+u^2x^3z+u^4xy,
\\[2ex]
T_2=&
-y^3zx^2-2y^2z^2x^2-yz^3x^2+x^3uy^2+2y^2x^2u^2-yz^3u^2-
\\[1ex] &
2u^2y^2z^2-y^3zu^2+xu^3z^2+2u^2x^2z^2+x^3uz^2+xu^3y^2-4x^2yzu^2
\\[1ex] &
-2x^3uyz+4xuy^2z^2+2xy^3zu-2xyu^3z+2xyz^3u.
\end{array}
$$
Из (\ref{facR}) видно, что если $C_3=0$ и $T\neq 0$, то
$Cycl=0$. К сожалению, утверждение $T\neq0$ неверно
(см.\ ниже). Система, рассматриваемая в
следующей Лемме, не эквивалентна (\ref{fsys}), но
будет ис\-поль\-зо\-ва\-на в \S~6
при доказательстве несуществования примыкающих к границе
се\-мейств решений системы (\ref{fsys}).

\begin{Lemma}
Всякое решение системы {\rm(\ref{fsys})}
на множестве $\Quad$
является также решением системы
\begin{equation}
\label{ZT}
T=0, \qquad T_z=0.
\end{equation}
\end{Lemma}

\Proof $\;$
Имеем
$\,T\cdot Cycl$=$F \cdot\Xi$,
где $\Xi$ --- рациональная функция без особенностей на $\Quad$.
Из уравнений (\ref{fsys}) следует
$$
 T\cdot Cycl\;=\;\partial_z (T\cdot Cycl)\;=0.
$$
Предположим сначала, что  $Cycl=0$, $T\neq 0$. Тогда
$\,\partial_z Cycl=0$. Но $\,\partial_z Cycl=-y\neq0$, противоречие.

В случае $\,Cycl\neq 0$, $T=0\,$ имеем $\,T_z=0$, как утверждалось.
Остается доказать, что невозможен случай
$T=Cycl=0$.
В результате подстановки $u=1$, $x=yz$ в $T$ получим
факторизующийся многочлен
$$
4\,y(y-1)(y+1)(y^2+2ty+1)\,z(z-1)(z+1)(z^2-2tz+1),
$$
который не обращается в 0 на множестве $\Quad$.
\BlackBox

\medskip\noindent
{\bf Следствие.}$\;${\it
Всякое решение системы {\rm(\ref{gradsys})}
на множестве $\Quad$
является также решением системы
}
\begin{equation}
\label{gradT}
T\,=\,T_x\,=\,T_y\,=\,T_z\,=\,T_u\,=\,0.
\end{equation}
Можно было бы надеяться доказать отсутствие решений системы
(\ref{gradsys}), доказав отсутствие решений полиномиальной системы
(\ref{gradT}). Но последняя имеет решения. (Например,
(\ref{gradT}) обращается в тождество при $u=x$, $y=z$.)
$\;$ Доказательство в
следующем п.\ основано на анализе
аналогичной системы с частичной рационализацией,
\begin{equation}
\label{gradR}
R\,=\,R_x\,=\,R_y\,=\,R_z\,=\,R_u\,=\,0,
\end{equation}
которая эквивалентна (\ref{gradsys}) на $\Quad$,
поскольку $\;F=R/F^*$, $\;F^*\neq 0$.

\subsection{Несуществование решений системы (\ref{gradR})}

\newcommand{\CR}{R}
\newcommand{\Rx}{\hat{R}_x}
\newcommand{\Ry}{\hat{R}_y}
\newcommand{\Rz}{\hat{R}_z}
\newcommand{\Ru}{\hat{R}_u}
\newcommand{\Rxx}{{Res}_x}
\newcommand{\Ryy}{{Res}_y}
\newcommand{\Rzz}{{Res}_z}
\newcommand{\Ruu}{{Res}_u}
\newcommand{\Rfx}{\tilde \Rxx}
\newcommand{\Rfy}{\tilde \Ryy}
\newcommand{\Rfz}{\tilde \Rxx}
\newcommand{\Rfu}{\tilde \Ruu}
\newcommand{\Rxy}{{Res}_{xy}}
\newcommand{\Ryu}{{Res}_{yu}}
\newcommand{\Ruz}{{Res}_{uz}}
\newcommand{\Rzx}{{Res}_{zx}}

Прежде всего, запишем производные $R_x$ и т.д.\ в
виде рациональных функ\-ций по $x,\dots,t$ и $\lambda$.
Имеем
$$
\lambda_x=\lambda\left(
\frac{a_x}{a}\,+
\frac{d_x}{d}\right)
=\lambda\left(\frac{x-yt}{a^2}\,+\,\frac{x+zt}{d^2}\right).
$$
Поэтому
$$
\Rx\;\eqbydef\;a^2 d^2\,\partial_x\CR
$$
--- алгебраический полином по переменным $x,y,z,t,\lambda$.
Аналогично, \mnog{}ами (по тем же переменным) являются
$$
\Ry\,=\,a^2 b^2\,\partial_y\CR,
\qquad
\Rz\,=\,c^2 d^2\,\partial_z\CR.
\qquad
\Ru\,=\,b^2 c^2\,\partial_u\CR,
$$
Получена система 5 полиномиальных уравнений
с неизвестными $x$, $y$, $z$, $u$, $t$ и $\lambda$
\begin{equation}
\label{gradRL}
\CR=\Rx=\Ry=\Rz=\Ru=0.
\end{equation}
Игнорируя уравнение связи (второе уравнение в (\ref{ratR})),
будем считать $\lambda$  независимой переменной от $\,x,\dots,t$.
Этим мы разрушаем однородность функции $\CR$ по переменным
$x,\dots,u$, следовательно 5 уравнений системы (\ref{gradRL})
более не являются зависимыми в силу тождества Эйлера.
Система (\ref{gradRL}) квазиоднородна по 5 переменным (исключая
$t$):  переменные $x,\dots,u$ имеют вес 1, а вес $\,\lambda$ равен
4. Одно из значений можно выбрать произ\-воль\-но.
Зафиксируем значение $u=1$. Теперь (\ref{gradRL})
становится системой 5 ура\-вне\-ний с 5 неизвестными $x$, $y$, $z$, $t$,
$\lambda$.

Доказать отсутствие решений этой небольшой полиномиальной си\-сте\-мы
оказалось нелегко. Автору пришлось перепробовать несколько схем
ис\-клю\-че\-ния.  Действуя "как попало", приходим к исчерпанию
вы\-чис\-ли\-те\-ль\-ных ресурсов.

Описываемая ниже схема оказалась выполнимой и, следовательно,
от\-но\-си\-тель\-но экономной, но все же далеко выходит за пределы
воз\-мож\-нос\-тей руч\-ного счета.

\medskip\noindent
\underline{Шаг 1}.$\;$
Исключаем $\lambda$, вычисляя результанты линейных
по $\lambda$ уравнений
$$
\Rxx\;\eqbydef\; \Resl{\CR}{\Rx}{\lambda},
$$
и аналогично $\Ryy$, $\Rzz$, $\Ruu$.
Результанты факторизуются, например,
$$
\Rxx=-4\,(t^2-1)\,(y+z)^2\,(x+1)\,\Rfx,
$$
где $\Rfx$ --- \mnog\ от $x$, $y$, $z$, $t$ с 92 членами.

\medskip\noindent
{\bf Замечание.}$\,$ На последующих шагах нам многократно будут
встречаться факторизации, в которых некоторые множители не
обращаются в 0 на множестве $Q$.
Например, таковы три множителя в формуле для $\Rxx$ выше;
другие примеры: $\,x$, $\;t\pm 1$, $\;z^2-2tzy+y^2$.
Такие множители будем называть тривиальными.
Некоторые множители, например, $\,y-1$, $x+1-y-z$, $y-xz$,  не
относятся к тривиальным, но приводят к простым ответвлениям от
основной линии. Мы рассматриваем соответствующие варианты в
леммах после завершения наиболее трудной части до\-ка\-за\-тель\-ства.
В тексте мы называем такие множители простыми (по английски
--- "simple", но не "prime").

\medskip\noindent
\underline{Шаг 2}.$\;$
Исключаем $t$. Факторизация результанта \mnog{}ов
$\Rfx$ и $\Rfy$
содержит, помимо тривиальных множителей, простые множители
\begin{equation}
\label{simp1}
z-xy,\quad y-xz,\quad x-yz,\quad x-y-z+1,
\end{equation}
и один "большой" множитель степени 6 с 18 членами, который обозначим
$\,\Rxy(x,y,z)$.
Аналогично, факторизация результанта $\Rfy$ и $\Rfu$ от\-но\-си\-тель\-но
$t$ приводит к тем же простым множителям и большому мно\-жи\-те\-лю
$\,\Ryu(x,y,z)$. Между коэффициентами многочленов
$\,\Rxy\,$ и $\,\Ryu\,$ имеется очевидное соответствие,
отвечающее симметрии исходной задачи относительно  пе\-ре\-ста\-нов\-ки
$\,x\rightarrow y\rightarrow u\rightarrow z\rightarrow x$,
$\,t\rightarrow -t$.

Аналогично определяются большие множители
$\,\Ruz\,$ и $\,\Rzx\,$ ре\-зуль\-тан\-тов
(относительно $t$) пар $\,(\Rfu,\,\Rfz)\,$ и $\,(\Rfz,\,\Rfx)$.

\medskip\noindent
\underline{Шаг 3} (решающий).
Имеют место факторизации (найдены методом проб с помощью \CAS).
$$
\begin{array}{c}
z\,\Rxy\,+\,y\,\Rzx\,=\,(y-z)(y+z)\,P^{(1)}_{4}(x,y,z), \\[2ex]
y\,\Ruz\,-\,z\,\Ryu\,=\,x(y-z)(y+z)\,P^{(2)}_{4}(x,y,z), \\[2ex]
P^{(1)}_{4}+P^{(2)}_{4}\,=\,(x-1)(x+1)(3x^2-2x+3+y^2+10yz+z^2).
\end{array}
$$
Последний множитель --- положительно определенная квадратичная фор\-ма.
Таким образом, остается рассмотреть "простые" случаи, когда
обра\-ща\-ет\-ся в $0$ либо один
из простых множителей (\ref{simp1}), либо $y-z$, либо $x-1$.


\medskip\noindent
\underline{Случай $x-y-z+1=0$}.$\,$   Имеем
$$
 \Rfy[z=x+1-y]=-4(y-1)(y^2+2yt+1)(x+1)^2(x-y)(x^2-2xyt+y^2).
$$
Следовательно либо
$y=1$, $x=z$, либо $x=y$, $z=1$.
В обоих случаях $Cycl=xu-yz=0$. Но равенство $F=F_z=0$ невозможно
на множестве $\Cycl$ (см.\ доказательство Леммы 5.2).

\pagebreak
\noindent
\underline{Случай $x=yz$}.$\,$
Аналогично, из факторизации
$$
 \Rfx[x=yz]=2(yz+1)\,y^2 (y^2-1)(y^2+2ty+1)\,z^2(z^2-1)(z^2-2tz+1)
$$
находим, что либо $y=1$, $x=z$, либо $z=1$, $x=y$, т.е.\
$Cycl=0$, и заключение как в предыдущем случае.

\medskip\noindent
\underline{Случай $y=xz$}$\,$ (Случай $y=xz$ симметричен.)
$\,$
Имеем факторизации
$$
\begin{array}{l}
 \Rfx[z=xy]=x^2(x^2-1)\,(y^2-1)(y^2+2ty+1)\,P_5(x,y,t),
 \\[1ex]
 \Rfu[z=xy]=-x(x^2-1)\,(y^2-1)(y^2+2ty+1)\,P_6(x,y,t),
 \end{array}
$$
где
$$
\begin{array}{l}
P_5(x,y,t)=x+x^2+y^2+y^4+x^2y^2-xy^4-4xy^3t,\\[1ex]
P_6(x,y,t)=1+x+y^2-xy^4+x^2y^2+x^2y^4-4xy^3 t.
 \end{array}
$$
Далее,
$$
 P_6-P_5=(x^2-1)(y^2-1)(y^2+1).
$$
Помимо подслучаев, ведущих к равенству $Cycl=0$ немедленно,
остается только случай $x=1$, $y=z$. При этом $a=c$, $b=d$,
и \chet{} $ABCD$ --- параллелограм. Выражение $C_3$ упрощается:
$\;C_3=2(abp-abq)$. Получаем $p=q$, откуда $x=y=z=u$, т.е.\
снова $Cycl=0$.

\medskip
Этим заканчивается разбор всех случаев и доказательство
не\-су\-ще\-ство\-ва\-ния решения системы (\ref{gradR}), а, значит,
и системы (\ref{gradsys}). 

\section{Несуществование решений системы (\ref{fsys})
вблизи границы множества $\Quad_1$}

\subsection{Стратификация границы}
Гипотетически, множество $\Mr$ (см.\ с.~11) может иметь предельные точки
на границе множества $\Quad_1$. Рассмотрим стратификацию границы
(раз\-би\-е\-ние ее на множества разных размерностей, открытые в
соответствующих ко\-ор\-ди\-нат\-ных подпространствах).
Грани коразмерности $k=1,\dots,4$ со\-от\-вет\-ству\-ют обращению ровно
$k$ из неравенств (\ref{rest5}) в равенства. (Все 5 не\-ра\-венств не могут
одновременно обратиться в равенства в силу усло\-вия (\ref{Quad1}).)

Введем терминологию, которая поможет нам компактно
описать мно\-го\-численные случаи. Пусть $\mathcal{V}$ --- некоторое
подмножество 5-элементного мно\-же\-ства переменных $x,\,y,\,z,\,u,\,t$.
 Будем говорить, что имеет место {\it случай $\mathcal{V}$},
если множество $\Mr$  имеет предельную точку
на грани, где неравенства (\ref{rest5}) для переменных из
$\mathcal{V}$ обращаются в равенства.

Например, имеет место случай $(x)$, если у системы (\ref{fsys}) есть
семейство решений, для которого $x\to 0$, а $y$, $z$, $u$ не
стремятся к 0 и 
$t\not\to\pm 1$.

\medskip
Перечислим случаи, которые нужно рассмотреть и исключить,
чтобы доказать Теорему 2.

\medskip\noindent
1) Случаи коразмерности 1: $\;$ $(x)$, $(y)$, $(z)\,$ и $\,(u)$.
Случай $\,(z)\,$ уже исключен Предложением 4.2. Остальные случаи
невозможны по симметрии, т.к.
\linebreak
урав\-не\-ние $\,F_z=0\,$ не
используется в Предложении 4.2.

\smallskip\noindent
2) Случай $(t)$ коразмерности 1. Он исключается Предложением
6.1.

\smallskip\noindent
3) Случаи $\,(xu)$, $(zy)\,$ и их вырождения
$\,(xut)$, $\,(zyt)\,$ покрываются Леммой 6.2.

\smallskip\noindent
4) Случай $(xy)$ исключается Предложением 6.3. По симметрии,
также исключается  случай $(uy)$.

\smallskip\noindent
5) Случай $(xz)$ и симметричный $(uz)$
исключаются Предложением 6.4.

\smallskip\noindent
6) Случаи коразмерности 2, когда одна
 из переменных --- $t$, например, $(xt)$, --- охватываются
 Пред\-ло\-же\-ни\-ем~6.5.

\smallskip\noindent
7) Случаи коразмерности 3 с тремя длинами: $(xyz)$, $(yzu)$, $(zux)$
и $(uxy)$ --- исключаются Леммой~6.6.
Оказываются охваченными и их вырождения коразмерности 4: $(xyzt)$ и
т.д.

\smallskip\noindent
8) $(xyt)$ и симметричный случай $(uyt)$.\\
9) $(xzt)$ и симметричный случай $(uzt)$.
Последние два случая оставлены за рамками настоящей работы.
Их анализ представляет самостоятельный интерес с точки зрения техники,
изложенной в  \cite[гл.~2]{Bru}. В этих случаях разрешение
особенности требует нескольких итераций.  Здесь мы
огра\-ни\-чи\-ва\-ем\-ся
формулировкой нужных результатов (Предложение 6.7).

\subsection{Случай $(t)$ (складывающийся \chet)}

\begin{Predl}
$\;$
 Случай $(t)$ невозможен.
\end{Predl}

\Proof $\;$
Без ограничения общности можно считать, что $\,t\to 1^-$.
Из (\ref{C3t1}) следует, что
равенство $\;\lim_{t\to 1} F =0\;$ невозможно, если \linebreak
$\lim (x-y)(u-z)>0$. $\;$ Имеется в виду предел вдоль некоторой кривой
$\,\langle x(t),y(t),z(t),u(t),t \rangle$,  $\;\;t\to 1^-$.
В дальнейшем из оставшихся двух воз\-можных комбинаций
знаков зафиксируем такую:
\begin{equation}
\label{limdiff}
\lim(y-x)\geq 0,\qquad
\lim(u-z)\geq 0.
\end{equation}

\noindent
Противоположный выбор ведет к полностью аналогичным вычислениям.

\noindent
Рассмотрим теперь 2 случая:

А)$\,$ Пределы (\ref{limdiff}) строго положительны;

Б)$\,$ Хотя бы один из этих пределов равен 0.

\medskip\noindent
\underline{Случай А}.$\,$
Если пределы (\ref{limdiff}) положительны, то с точностью
до членов порядка $\tau$ имеем
\begin{equation}
\label{tasym}
\begin{array}{lcl}
\dst a\sim (y-x)+\frac{xy}{y-x}\tau,&\quad& \dst
 b\sim (u+y)-\frac{uy}{u+y}\tau,\\[2ex]
\dst c\sim (u-z)+\frac{uz}{u-z}\tau, &\quad& \dst
 d\sim (z+x)-\frac{zx}{z+x}\tau.
\end{array}
\end{equation}
Отсюда находим асимптотику при $\tau\to 0$:
$$
 F\cdot abcd = \tau\,(y+z)(x+y)(xu-yz) \,P_3(x,y,z,u)  \,+\, O(\tau^2),
$$
где
$$
P_3=ux^2+yz^2+3yzx+3zux-xu^2-zy^2-3uyz-3uyx.
$$
Многочлен $P_3$ имеет степень 2 по переменной $x$. Коэффициент при
$x^2$ равен $u>0$,
и $\;\;P_3[x=0]\,=yz(z-y-3u)<0$,  $\;\;P_3[x=y]\,=\,y(z-u)\,(u+2y+z)<0$.
Следовательно, равенство $P_3=0$ несовместно с неравенствами
(\ref{limdiff}).  Остается возможность обращения $F$ в 0 при
$\tau>0$ в случае $xu-yz\to 0$.  Эта возможность реализуется ---
ср.\ условие (\ref{C2}).  Покажем, однако, что при этом $F_z\neq 0$.

Подставляя (\ref{tasym}) в $F_z$, находим асимптотику
при $\tau\to 0$:
$$
 F_z\cdot(a b c d) \sim -\tau\,(x+u)\, P_7(x,y,z,u),
$$
Многочлен $\,P_7\,$ однороден; подстановка $\,u=1\,$ и $\,x=yz\,$ приводит к
много\-члену, который факторизуется:
$$
 P_7(yz,y,z,1)=\,4 \,[y (y+1) z(z-1)]^2 \,(y+z).
$$
В рассматриваемом случае $\lim(z-1)=\lim(z-u)<0$. Следовательно,
$F_z\neq 0$ при малых $\tau>0$.

\medskip\noindent
\underline{Случай Б}.$\,$
Предположим, что
$ z-u\to 0$.
Переключаясь на рассмотрение многочлена $T$ (см~.п.~5.2) вместо
$F$, имеем простую факторизацию:
\begin{equation}
\label{faczequ}
T[z=u,\, t=1] \;=\,u\,(x-y)^3\,(x+y+2u).
\end{equation}
Следовательно, необходимо также $x-y\to 0$.
(Можно было пред\-по\-ло\-жить последнее условие и вывести первое).
Итак, имеем 3 малых параметра:
$$
\tau=1-t, \qquad
v=y-x, \qquad s=u-z
$$
 Заметим, что теперь мы не имеем права считать,
что $v$ и $s$ положительны.

Этот случай сильно вырожден и труден. Для его анализа мы будем
использовать как пару уравнений $F=F_z=0$, так и пару уравнений
$T=T_z=0\;$ --- см.\ Лемму 5.2.

\medskip\noindent
$1^\circ$. Обратимся сначала к величине $F=C_3$. Запишем
\begin{equation}
\label{CAC}
 C_3=Aa+Cc,
 \qquad A=bp-dq,
 \qquad C=dp-bq,
\end{equation}
и вычислим асимптотику коэффициентов $A$ и $C$.
В данном случае
\begin{equation}
\label{tasbd}
 b=(u+y) -\frac{uy}{u+y}\,\tau + O(\tau^2),\qquad
 d=(x+z) -\frac{xz}{x+z}\,\tau + O(\tau^2).
\end{equation}
Обозначая
$$
  \xi=|s|+|v|
$$
и подставляя асимптотики (\ref{tasbd}) в выражения для $A$ и $C$,
находим
\begin{equation}
\label{tasAC}
 A= 2(x+u)s + O(\tau \xi)+O(\tau^2),
 \qquad
 C=-2(x+u)v + O(\tau \xi)+O(\tau^2).
\end{equation}
В действительности члены, содержащие некоторую степень $\,\tau\,$
без 
мно\-жи\-те\-лей $s$ или $v$,
отсутствуют, но мы докажем этот факт не для величин
$A$ и $C$ в отдельности, а для уравнения $\,Aa+Cc=0$.
До\-ка\-за\-тель\-ство просто: при $\,v=s=0\,$ имеем $\,C_3=0\,$
(геометрически --- это случай равно\-боч\-ной трапеции).

Возводя уравнение $\,Aa=-Cc\,$ в квадрат и используя
формулы
\begin{equation}
\label{tasac}
 a^2= v^2+2x y\tau,
 c^2= s^2+2u z\tau,
\end{equation}
получаем после сокращений
$$
 x^2 s^2 - u^2 v^2 = O(\xi^2(\xi+\tau)).
$$
Следовательно, $s\sim \pm (u/x)v$.
Возвращаясь к 
уравнению $\,Aa=-Cc$, определяем, что ему
удовлетворяет только ветвь со знаком $+$.
(Полу\-чен\-ный результат можно было ожидать --- он согласуется
с условием (\ref{Cycl0}).)

\medskip\noindent
$2^\circ$.
Теперь обратимся к уравнению $\partial_zC_3 =0$.
Используя соотношения $\;c\,\partial_z c=z-ut\;$
и $\;d/dz=-d/ds$ и независимость $a$ от $z$, выводим из (\ref{CAC}):
$$
c\,\partial_zC_3 \;=\; X+ ac Y,
$$
где
$$
 X=-c^2\,\partial_s C+ (z-ut)C,
 \qquad
 Y=\partial_s A.
$$
Из (\ref{tasAC}) видно, что $\lim Y =2(x+u)>0$,
а из (\ref{tasac}) следует, что
\begin{equation}
\label{ineqactau}
ac>\,\const\, \tau .
\end{equation}
Следовательно,
$\;\tau=O(ac\,Y)$.  
Обратимся к коэффициенту $X$.
Поскольку $\,\partial_s C=O(\xi)$, первое слагаемое в $X$ есть
$\,O(\xi^3+\xi\tau)$. Второе слагаемое
$$
(z-ut)C = (-s+u\tau)C = 2(x+u)vs + O(\xi^3+\xi\tau).
$$
Таким образом, на асимптотическом решении уравнения
$\,X+ ac Y=0\,$  должна выполняться оценка
$$
\tau=O(X)=O(\xi^2).
$$
Эту оценку можно усилить, поскольку неравенство
(\ref{ineqactau}) слишком грубое.
В силу результата п.$1^\circ$, $\;\xi^2$ и $vs$ имеют одинаковый
асимптотический порядок.
Если предположить, что $\;\lim \tau/\xi^2 >0$, то получим
$$
\lim \frac{ac\,Y}{vs} > 2(x+u) = \lim \frac{X}{vs},
$$
противоречие. Следовательно,
$\,\tau=o(\xi^2)$.

\medskip\noindent
$3^\circ$.
Укорочение уравнения $\;T=0\;$ при малых $\,s$, $v$, $\tau\,$ имеет вид
$$
 T\sim 4(x+u)^2\,(s-v)(s+v)(sx+v) + O(\tau\xi+\xi^4).
$$
Первое слагаемое доминирует, и мы знаем, что $sx/vu\to 1$.
Следовательно, на асимптотическом решении $\;x/u\to 1$.

\medskip\noindent
$4^\circ$.
Обратимся, наконец, к уравнению $T_z=0$ и найдем его укорочение
при $\,s,\,v,\,\to 0\,$, $\,w=x-u\to 0\,$ и $\tau=o(s^2)$:
$$
 T_z\,\sim\, -16 (s+v)(3s-v)+ o(\xi^2).
$$
Равенство $T_z=0$ несовместно с выведенным в
$3^\circ$ условием $s/v\to 1$.
Полученное противоречие доказывает невозможность Случая Б.
\BlackBox

\subsection{Оставшиеся предельные случаи}

\begin{Lemma}
{\rm(О сплющивающемся четырехугольнике).}$\,$
 Случаи $(xu)$, $(zy)$ и их вырождения $(xut)$, $(zyt)$
 невозможны.
  Более того, одно лишь урав\-не\-ние $F=0$ не имеет семейств решений,
  выходящих на со\-от\-вет\-ству\-ю\-щие грани.
\end{Lemma}

\Proof $\;$
Предположим, что $x\to0$, $u\to 0$, а $z$ и $y$ не малы.
Тогда $\,a,b\to y$, $\,c,d\to z$, $\,p\to 0$.
Следовательно, $\,F\to -2yz(y+z)\neq 0$.
\BlackBox
\begin{Predl}
{\rm (Вырождение четырехугольника в треугольник -- I)}$\;$
Случай $(xy)$ невозможен.
\end{Predl}

\Proof  $\;$
При $x\to 0$, $y\to 0$ имеем
$$
a=O(x+y),
\;\quad
b=u+yt+O(y^2),
\qquad
c=\const,
\;\quad
d=z+xt+O(x^2).
$$
Отсюда находим асимптотики
\begin{equation}
\label{sisxy}
\begin{array}{c}
F = x\alpha+ y\beta + a\gamma\,+\,O(x^2+y^2),
\\[1.8ex]
F_z = x\tilde\alpha+ y\tilde\beta + a\tilde\gamma\,+\,O(x^2+y^2),
\end{array}
\end{equation}
с коэффициентами
\begin{equation}
\label{coefsisxy}
\begin{array}{lccl}
\alpha=c(ut+z),&&
\tilde\alpha=z(z^2-u^2 t^2+c^2),
\\[1.8ex]
\beta=-c(u+zt), &&
\tilde\beta=(1-t)z^3+(1+t)(u-z)uz,
\\[1.8ex]
\gamma=u^2-z^2,&&
\tilde\gamma=-2 z^2 c.
\end{array}
\end{equation}
Найдем предельный коэффициент пропорциональности $x$ и $y$,
исключая
$a$ из системы (\ref{sisxy}). Получаем
$$
 \frac{x}{y}\,\sim\,\frac{\tilde\beta \gamma - \beta\tilde\gamma}%
 {\alpha\tilde\gamma-\tilde\alpha\gamma}
 \,=\,
 \frac{z}{u}\;\frac{3u^2+z^2}{3z^2+u^2}\,=:\,s_{*}.
$$
(В ходе выкладок использовали выражение $c^2$ через $u$, $z$,
$t$, см.~(\ref{elem6}).)

Завершим доказательство Предложения, показав, что
укорочение вто\-ро\-го из уравнений (\ref{sisxy}) не имеет решений
с $\, x/y=s_{*}$.
Избавляясь в этом уравнении от иррациональности
$\;a=(x^2+y^2-2xyt)^{1/2}$, получаем квад\-ра\-тич\-ную форму
$$
 (\tilde\alpha^2-\tilde\gamma^2)\,x^2\;+\;
 (\tilde\beta^2-\tilde\gamma^2)\,y^2\;+\;
 (2\tilde\alpha\tilde\beta+2 t\tilde\gamma^2) xy\;
 =\,0.
$$
Подставляя коэффициенты из (\ref{coefsisxy})
и полагая $\,u=1$, $y=1$, $x=s_{*}\,$ (выражения однородны!),
получим многочлен от $z$ с параметром $t$, до\-пус\-ка\-ю\-щий
факторизацию
$$ z^4\,(z^2-2tz+1)\,(3z^2-2tz+3).
$$
Он не имеет положительных корней, если $|t|<1$.
\BlackBox

\begin{Predl}
{\rm(Вырождение четырехугольника в треугольник -- II)}$\;$
Случай $(xz)$ невозможен.
\end{Predl}

\Proof  $\;$
Доказательство аналогично предыдущему. В конце надо будет
разобрать один более тонкий случай.

При $x\to 0$, $z\to 0$ имеем
$$
a=y-xt+O(x^2),
\qquad
b=\const,
\qquad
c=u-zt+O(z^2).
\qquad
d=O(x+z).
$$
Аналогично (\ref{sisxz}), находим
\begin{equation}
\label{sisxz}
\begin{array}{c}
F = x\alpha+ z\beta + d\gamma\,+\,O(x^2+z^2),
\\[2ex]
F_z = x\tilde\alpha+ z\tilde\beta + d\tilde\gamma\,+\,O(x^2+z^2),
\end{array}
\end{equation}
с коэффициентами
\begin{equation}
\label{coefsisxz}
\begin{array}{lccl}
\alpha=b(y-ut),&&
\tilde\alpha=ut(u^2 -y^2),
\\[1.8ex]
\beta=-b(u-yt), &&
\tilde\beta=u(u^2 -y^2),
\\[1.8ex]
\gamma=u^2-y^2, &&
\tilde\gamma=bu(ty-u).
\end{array}
\end{equation}
Предельный коэффициент пропорциональности $x$ и $z$ равен
$$
s_{*}=\,
 \frac{y}{u}\;\frac{3u^2-2tuy-y^2}{y^2-2tuy+u^2}.
$$
Рассматривая квадратичную форму, соответствующую
первому из урав\-не\-ний (\ref{sisxz}), получим после упрощений
и подстановок $\,u=1\,$,  $\,x=s_{*}z\,$ пре\-дель\-ное уравнение для
$z$ и $t$
$$
  4(t^2 - 1)\,y^2(y^2-1)^2\,(y^2+2ty-3)=0.
$$
Таким образом, остается возможность существования асимптотического
решения, на котором
\begin{equation}
\label{xzdegen}
\;y^2+2ty-3\to 0.
\end{equation}
Рассмотрение квадратичной формы, соответствующей
второму урав\-не\-нию в (\ref{sisxz}) также не исключает этот случай.
Заметим, что в этом слу\-чае $s_{*}=0$, т.е. $x=o(z)$. Следовательно,
$\; d=z+o(z)$. Получаем (с $u=1$)
\begin{equation}
\label{xzdegen2}
F = z\,((1-b)+bty-y^2) \,+\,o(z).
\end{equation}
Условие (\ref{xzdegen}) влечет $\,b\to 2$. Коэффициент при $z$ в
(\ref{xzdegen2}) в пределе равен
$\;-y^2+2ty-1 <0$.
Этим возможность (\ref{xzdegen}) исключается.
\BlackBox

\begin{Predl} %
Случаи $(xt)$, $(ut)$, $(yt)$, $(zt)$ невозможны.
\end{Predl}

\Proof  $\;$
Рассмотрим первый из двух подслучаев случая
$(xt)$: $\;x\to 0$, $\,t\to+1$.
Вначале, следуя доказательству Предложения 6.1,
по\-лу\-чим (принимая во внимание, что $y>x$)
\begin{equation}
\label{xt_inequz}
\lim (u-z)\geq 0.
\end{equation}
Рассмотрим уравнение $T=0$.
При $t=1$, $x=0$ оно факторизуется:
$$
 T[t=1,\;x=0]\,=\,-yz(3u+y-z)(u+y-z)(-u+y+z)(u+y+z)=0.
$$
Ввиду (\ref{xt_inequz}), факторы $\,(3u+y-z)\;$ и $\;(u+y-z)\,$
не могут в пределе обращаться в $0$, поскольку $y$ ограничено снизу.
 Остается единственная возможность
$$
\lim (y+z-u)=0.
$$
Однако
$$
 T_z[t=1,\;x=0,\;u=y+z]\,=\,-8y^2 z(y+z)(2y+z)\neq 0.
$$
Подслучай $t\to -1$ приводит к факторизации
$$
 T[t=-1,\;x=0]\,=\,-yz(-3u+y-z)(-u+y-z)(-u+y+z)(u+y+z),
$$
откуда снова следует, что на предельном семействе
$\;u-y\to z$. Но
$$
 T_z[t=-1,\;x=0,\;u=y+z]\,=\,-8y z^2(y+2z)(y+z)\neq 0.
$$
Случай $(ut)$ не требует отдельного доказательства ввиду симметрии
$x\leftrightarrow u$, $t\leftrightarrow 1-t$.
Остальные случаи, описанные в Предложении, доказываются аналогично.
Приведем соответствующие факторизации.

\smallskip\noindent
Случай $\;y\to 0$, $\;t\to 1$.
\begin{equation}
\begin{array}{l}
 T[y=0, \;t=1]=xu(u+x-z)(-u+x+z)(u+x+z)(-u+x+3z),
\\[1ex]
 T_z[y=0, \;t=1,\;z=x+u]=-8x^2 u(u+x)(2x+u)\neq 0.
 \end{array}
\end{equation}

\smallskip\noindent
Случай $\;z\to 0$, $\;t\to 1$.
\begin{equation}
\begin{array}{lr}
 \dst T[y=0, \;t=1]=xu(u+y-x)(x+u-y)(u+x+y)(-u+x+3y),&
\\[1ex]
 \dst T_z[y=0, \;t=1,\;z=x+u]=-8x^2 u(u+x)(2x+u)\neq 0.&
 \end{array}
\end{equation}
Предложение доказано. \BlackBox

\begin{Lemma}
 Первое укорочение \mnog{}а $T$ при $u=1$,
 $\,x,y,z\to 0$, есть
$$
 T(x,y,z,1,t)\sim x+3yz.
$$
\end{Lemma}
\Proof  $\;$Прямая проверка. \BlackBox

\medskip
\noindent{\bf Следствие.}$\;$ {\it Случай $\mathcal{V}$,
когда множество $\mathcal{V}$ содержит три из четырех параметров
$x$, $y$, $z$, $u$, невозможен.}

\begin{Predl}
Система {\rm(\ref{fsys})} не имеет семейств решений с асим\-пто\-ти\-ками
$t\to \pm 1$, $\xi,\eta\to 0$, где $\xi$, $\eta$ --- два из
четырех параметров $x$, $y$, $z$, $u$.
\end{Predl}

\noindent
Это предложение в настоящей работе не доказывается.

\subsection*{Приложение. 
Нули функции $C_3$ для \emph{невыпуклых}\ сим\-мет\-рич\-ных \chet{}ов с
пер\-пен\-ди\-кулярными диагоналями
}


Фиксируем масштаб: $z=1$.
Полная система ограничений есть
$$
 u=x,\qquad t=0,\qquad z=1.
 \eqno({\rm A}.1)
$$
При этих условиях
$$
 C_3=2x(a^2+d^2)-2(y+1)ad,
$$
и рационализованное выражение $\,R\,$ (\ref{ratR}) факторизуется
$$
\begin{array}{l}
\dst R\,=\,C_3\cdot(x(a^2+d^2)+(y+1)ad)
 \;=\;2(x^2-y)\,P_4(x,y),\\[2ex]
 P_4(x,y)=(3x^2+4x^4)+ (1+2x^2)y+ (2+3x^2)y^2+ y^3.
\end{array}
$$
При $-1<y<0$ рационализующий множитель положителен, поэтому
множества нулей функции $\,C_3|$(A.1)$\,$ и многочлена $\,P_4(x,z)\,$
при условии $-1<y<0$ совпадают.
Решение с максимальным значением $x$ находится из системы уравнений
$$
 P(x,z)=0, \qquad
 \partial_z P=(3z+1)(z-2x^2-1)=0.
$$
Решение с положительным $x$ единственно:
$$
 z=-1/3,\qquad
 x=\frac{1}{3}\sqrt{2\sqrt{3}-3} \approx 0.227083.
$$
При этом
$$
 a=b\approx 1.025459, \;\quad
 c=d\approx 0.403334, \;\quad
 p \approx 0.454167, \;\quad
 q \approx 0.66667.
$$




\begin{thebibliography}{2}
\bibitem{Berger} М.~Берже. Геометрия, т.1.
М.:~Мир, 1985.
%
\bibitem{Bru} А.Д.~Брюно. Степенная геометрия в алгебраических и
диф\-фе\-рен\-ци\-аль\-ных уравнениях.
М.:~Наука, 1998.
%
%
\bibitem{Coxeter} Г.С.М.~Коксетер, С.Л.~Грейтцер.
Новые встречи с геометрией.
М.:~Наука, 1978.
%
\bibitem{Prasol} В.В.~Прасолов. Задачи по планиметрии, ч.2.
М.:~Наука, 1986.
%
\bibitem{Sabitov} И.Х.~Сабитов. Объемы многогранников.
М.:~Изд.~МЦНМО, 2002.
%
%
%
%
%
\end{thebibliography}
\end{document}